\documentclass{article}
\usepackage[latin1]{inputenc}
\usepackage{a4wide}
\usepackage{amssymb}
\usepackage[arrow, matrix, curve]{xy}
\usepackage{amsmath}
\usepackage{amsthm}
\usepackage{bbm}
\usepackage{paralist}
\usepackage{mathtools}

\newtheoremstyle{style}
  {0.3cm}
  {0.3cm}
  {}
  {}
  {\normalfont\bfseries}
  {:}{ }
  {{\normalfont\bfseries \thmname{#1}\thmnumber{ #2}\thmnote{ (#3)}}}
\newtheoremstyle{stylesatz}
  {0.3cm}
  {0.3cm}
  {\itshape}
  {}
  {\normalfont\bfseries}
  {:}{ }
  {{\normalfont\bfseries \thmname{#1}\thmnumber{ #2}\thmnote{ (#3)}}}
\newtheoremstyle{stylebeweis}
  {0.3cm}
  {0.3cm}
  {}
  {}
  {\normalfont\itshape}
  {:}{ }
  {{\normalfont\itshape \thmname{#1}\thmnumber{ #2}\thmnote{ (#3)}}}
\theoremstyle{style}
\newtheorem{Rem}{Remark}
\newtheorem{Def}[Rem]{Definition}
\newtheorem{Not}[Rem]{Notation}
\newtheorem{App}[Rem]{Application}
\newtheorem{Exa}[Rem]{Example}
\theoremstyle{stylesatz}
\newtheorem{Thm}[Rem]{Theorem}
\newtheorem{Lem}[Rem]{Lemma}
\newtheorem{Kor}[Rem]{Corollary}
\newtheorem{Alg}[Rem]{Algorithm}
\newtheorem{Prp}[Rem]{Proposition}
\theoremstyle{stylebeweis}
\newtheorem*{Bew}{Proof}
\numberwithin{Rem}{section}
\numberwithin{equation}{section}

\DeclarePairedDelimiter\abs{\lvert}{\rvert}
\DeclarePairedDelimiter\norm{\lVert}{\rVert}

\begin{document}

\title{The self-justifying Elo rating system}
\author{Fabian Langholf\\\bigskip\small http://fabian.langholf.de}
\date{}
\maketitle

\begin{abstract}
We suggest an improvement of the Elo rating system. Whereas Elo's theoretical background remains unaffected, we significantly change the way in which rating values are adjusted. It turns out that the modified system behaves much more naturally, and that it offers several advantages over the classical one. The key idea is a fixed point approach to the definition of the rating values. We provide an algorithm for the purpose of their computation.
\end{abstract}

\section{Introduction}
Elo rating systems -- as described in \cite{ELO} -- find widespread use in settings where the `strength' of some `players' is supposed to be estimated on the basis of pairwise competition results. Elo is mainly interested in chess, but many other sports (such as go and football) and several other disciplines apply his method with great success. The core of such a system is a mechanism that converts rating differences into performance expectations. Whenever the actual results do not meet the expectations, the rating is adjusted.

Elo develops a comprehensive theory of the `correct' mapping between rating and results. His approach is based on probability theory, he spends great effort on the choice of the best fitting probabilistic model and discusses aspects such as reliability and integrity (\cite{ELO}, e.g. sections 1.3, 2.5, 3.5, chapter 8). Furthermore, he shows in detail how the rating can be applied in chess (e.g. historical or geographical comparisons, lifetime development, titles). None of his results are contradicted in this paper.

What we address is of technical nature, namely the set of assignments Elo uses to define his precise rating values. His main formula (\cite{ELO}, section 1.6, formula (2)) is used to adjust the rating when new results are available. Realising that initial ratings and tournament performance ratings cannot be treated in the same way, Elo introduces an additional procedure with several variants (\cite{ELO}, section 1.5, formula (1)). He accepts lots of approximations and heuristics, certainly due to restricted computational power in his day. This leads to a series of peculiarities (some of which are discussed in the list below in this introduction).

One particular aspect is worth closer inspection. As mentioned above, rating adjustments are based on performance expectations. Elo's classical approach (see application \ref{Appcl} below) calculates performance expectations which are based on past results. The expectations are compared with the present results. This means that past and present results influence the new rating in completely different ways. What the classical Elo does not do at all is to check whether the resulting rating is coherent with the (past or present) competition results. So it can happen that the very same performance, repeated in two sequential tournaments, will lead to an increasing rating in the first case and a cancellation of the just gained points in the second case (example \ref{Exadiv}).

The approach discussed here (application \ref{Appsj}) overcomes the unequal treatment of past and present results. It provides a rating with perfect coherence between results and expectations: we adjust the rating by comparing the competition results with the expectations based on the final, the already adjusted rating values (which might sound like magic at first). Since our rating adjustments are `justified' by the expectations based on their own resulting rating values, we say that the system is self-justifying.

We will see that a unique self-justifying rating always exists (theorem \ref{hauptthm}), and that it inherits basic properties of the classical Elo rating, such as continuity and the maximal amplitude of rating adjustments (proposition \ref{Prpcont}). Both rating approaches behave similarly when they are confronted with competition results of small magnitude (or, equivalently, with small values of the parameter $k$ to be introduced in the main text; proposition \ref{Prpsmall}). More importantly, in case of convergence of the classical Elo, the self-justifying Elo converges to the same rating values (corollary \ref{Korconv}). In section \ref{SctCalc}, we show how the self-justifying rating can be computed.

In the following list, we summarise the main advantages of the self-justifying Elo over the classical one.
\begin{enumerate}[(1)]
\item{The self-justifying rating does not depend on previously computed rating values (but only on the competition results). It does neither depend on the order of the results nor on the update cycles.}
\item{If the set of players can be divided into two subsets, where the members of one subset always defeat those of the other one, we have to expect unbounded growth of their rating values (in either rating system). In all other cases, repetition of competition results leads to convergent self-justifying Elo ratings (theorem \ref{Thmbigk}) -- in contrast to the classical situation (example \ref{Exadiv}).}
\item{If a player defeats another one, then this will lead to a better self-justifying rating (proposition \ref{Prpmon}) -- the opposite can be true in the classical case (example \ref{Examon}).}
\item{The classical rating tends to overcompensation when adjusting rating differences (example \ref{Exaoco}). To avoid this, the impact of a single competition result on the rating must be minimized (i.e. relatively small values of the parameter $k$ to be defined must be chosen). The self-justifying rating behaves better, so no such restrictions must be accepted.}\label{advhk}
\item{This allows us to use the self-justifying Elo to rate single tournaments (`performance ratings' in the sense of \cite{ELO}, section 1.5). A differing, heuristical procedure as in the classical case can therefore be avoided.}
\item{Players may be tempted to choose their opponents selectively, avoiding underrated and preferring overrated ones. The self-justifying rating, however, is endowed with a mechanism of self-correction, making this behaviour less attractive (remark \ref{Remsel}).}
\item{Some implementations of the classical Elo (among these the original one in \cite{ELO}) include complex processes dealing with initial ratings of new players. Because of the self-correction mechanism and advantage (\ref{advhk}), the self-justifying Elo does not require anything like this.}
\item{Weights for the competition results can easily be defined in the self-justifying case (application \ref{Appsj}). It is much more subtle to understand how a specific result affects the classical rating -- we have already noted above that these `implicit weights' can even be negative (example \ref{Examon}).}
\item{If the rating values of retiring players systematically differ from the average (or when special mechanisms for initial ratings are implemented), then the average rating of active players will inflate or deflate. The self-justifying Elo offers two simple possibilities to prevent this, either by simply removing inactive players from the calculation or by using decreasing weights (application \ref{Appsj}).}
\end{enumerate}

It should be clarified that Elo's original rating system as a whole -- including all of its complex components -- masters most of the issues above. But the bottom line is that that there exists a more natural approach, making much of the complexity obsolete.

\section{Preliminaries}

\begin{Not}[Ratings, competition results]\label{notation}
We fix $n\in\mathbb{N}$ with $n\geq2$ (the `number of competing players') and make use of the identification $n=\{0,1,...,n-1\}$ (the `set of competing players'). Ratings are modelled as elements of the hyperplane $\mathcal{R}:=\{x\in\mathbb{R}^n\ |\ \sum_{i\in n}{x_i}=0\}$ (the `rating space', endowed with the norm $\norm*{.}_1$). We introduce $\mathcal{P}:=\{p\in\mathbb{R}^{n\times n}\ |\ \forall i,j\in n:p_{ij}\geq0,\ \forall i\in n:p_{ii}=0\}$ (the `space of competition results', endowed with the metric induced by the norm $\norm*{.}_1$ on $\mathbb{R}^{n\times n}$; $p_{ij}$ is interpreted as the number of points that player $i$ gained against player $j$).
\end{Not}

\begin{Def}[Classical Elo map]\label{defclelomap}
Let $k>0$. The \emph{classical Elo map} with dynamising parameter $k$ is
\begin{equation*}
elo_k^{cl}:\mathcal{R}\times\mathcal{P}\rightarrow\mathcal{R},\ (x,p)\mapsto k\cdot\left(\sum\limits_{j\in n}{\left(p_{ij}-\frac{p_{ij}+p_{ji}}{1+exp(x_j-x_i)}\right)}\right)_{i\in n}
\text{.}
\end{equation*}
\end{Def}

\begin{Rem}[Symmetry]\label{Remsym}
Let $x\in\mathcal{R}$, $p\in\mathcal{P}$, $k>0$, $i,j\in n$. Using
\begin{equation}\label{expconst}
\frac{1}{1+exp(t)}+\frac{1}{1+exp(-t)}=\frac{1}{1+exp(t)}+\frac{1}{1+\frac{1}{exp(t)}}=\frac{1}{1+exp(t)}+\frac{exp(t)}{1+exp(t)}=1
\end{equation}
for $t\in\mathbb{R}$, we see that the $j$-th summand of $\left(elo_k^{cl}(x,p)\right)_i$ and the $i$-th summand of $\left(elo_k^{cl}(x,p)\right)_j$ differ only by sign. This implies that $\sum_{i\in n}{\left(elo_k^{cl}(x,p)\right)_{i}}=0$, showing that the classical Elo map is well-defined. Furthermore, our observation contributes to the interpretation of the map. If player $i$ and player $j$ compete for one point, then the term $\frac{1}{1+exp(x_j-x_i)}$ is meant to be the number of points that player $i$ can expect, reflecting the difference of the ratings of both players. The Elo map compares the gained points with the expected points.

Obviously, the map $\mathcal{R}\times\mathcal{P}\times\mathbb{R}_{>0}\rightarrow\mathcal{R},\ (x,p,k)\mapsto elo_k(x,p)$ is continuous. Furthermore,
\begin{equation}\label{absch_2kp}
\norm*{elo_k^{cl}(x,p)}_1
\leq k\cdot\sum\limits_{\substack{i,j\in n\\i<j}}{2\cdot\abs*{p_{ij}-\frac{p_{ij}+p_{ji}}{1+exp(x_j-x_i)}}}
\leq2\cdot k\cdot\sum\limits_{\substack{i,j\in n\\i<j}}{\left(p_{ij}+p_{ji}\right)}
=2\cdot k\cdot\norm*{p}_1
\text{.}
\end{equation}
\end{Rem}

\begin{Rem}[Embeddings]\label{Rememb}
By notation \ref{notation}, the number $n$ of competing players is fixed. However, it is clear that rating spaces and spaces of competition results defined on subsets of $n$ canonically embed into the ones defined above, and that the definition of Elo is compatible with these embeddings.
\end{Rem}

\begin{App}[Classical Elo]\label{Appcl}
We discuss the classical way to implement Elo rating systems. We choose $k>0$, $\mu\in\mathbb{R}$ (the `average rating'), and $\sigma>0$ (the `deviation factor'). We set $x^0:=(0)_{i\in n}\in\mathcal{R}$. We decide how often to update the rating (e.g. monthly or immediately after a new result is known). This decision leads to a (finite) sequence $p^0,p^1,...,p^{m-1}\in\mathcal{P}$ of competition results, corresponding to consecutive time periods. For $l\in m$, the classical Elo rating at the end of period $l$ is
\begin{equation*}
x^{l+1}:=x^l+elo_k^{cl}(x^l,p^l)\in\mathcal{R}\text{.}
\end{equation*}
In particuar, $x^m$ is the current classical Elo rating.

Since the specifications of the rating space and the Elo map concerning average and deviation might not match the requirements / the conventions, an affine transformation is performed at the end. The result is $(\mu)_{i\in n}+\sigma\cdot x^m$ (the `published rating').
\end{App}

\begin{Rem}[Parameter $k$]\label{Remsuppk}
The parameter $k$ determines to what extent a competition result changes the rating values. A small value of $k$ will lead to a rigid rating system, unable to induce the correct rating differences between the players. If a large value is used, then excessive adjustments are performed, undermining the validity of the rating.

Obviously, if we replaced all competition results $p$ by $k\cdot p$, we could suppress the parameter $k$ from the notation. It turns out, however, that several results can be expressed most conveniently using the parameter. Furthermore, most practical applications are associated with natural ideas of what a point is, and the parameter allows us to stick to these natural scales.
\end{Rem}

\begin{Rem}[Original presentation]
Elo would complain about our presentation above, because we seemingly ignore several aspects that are important to him. We give the link between our `classical Elo' and the original system in this remark.

The essential basis for our definitions is \cite{ELO}, section 1.6, formula (2). Elo does not vary our parameter $\sigma$, but uses the fixed value $\sigma=\frac{400}{\ln{(10)}}$. He uses the expression `class interval' for a rating difference of $\frac{\ln{(10)}}{2}$ (corresponding to $200$ points on his scale), and finds a natural meaning for this notion (\cite{ELO}, sections 1.2, 1.3). In contrast to our description, he does not use a fixed average value (such as our $\mu$): the differing procedure to provide initial ratings and the actions to prevent deflation (\cite{ELO}, chapter 3) are not compatible with a constant average. Since the self-justifying Elo does not need to cover these aspects, we do not focus on them.

Our parameter $k$ corresponds to Elo's parameter $K=\sigma\cdot k$. Elo uses varying values of $K$ in the interval between $K=10$ and $K=32$, depending on the `player development' (\cite{ELO}, 1.63, 8.28).

Our definition \ref{defclelomap} uses the logistic distribution function (\cite{ELO}, section 8.4). Elo discusses other possibilities. In fact, another function (normal distribution, \cite{ELO}, section 1.3) is his first choice, but he already recognizes some advantages of the logistic function.

In practice, there exist various more variants of the procedure described in application \ref{Appcl}.
\end{Rem}

\section{Construction}

\begin{Def}[Self-justifying rating]
Let $k>0$ and $p\in\mathcal{P}$. A rating $x\in\mathcal{R}$ is said to be \emph{self-justifying} with respect to $p$ and $k$ if it is a fixed point of the classical Elo map $elo_k^{cl}(.,p):\mathcal{R}\rightarrow\mathcal{R}$, i.e.
\begin{equation*}
x=elo_k^{cl}(x,p)\text{.}
\end{equation*}
\end{Def}

\begin{Rem}[Contraction function]\label{Remphi}
Let $k>0$ and $p\in\mathcal{P}$. Let $0\leq\xi<1$. Then the map 
\begin{equation}\label{phidef}
\varphi_\xi:\mathcal{R}\rightarrow\mathcal{R},\ x\mapsto\xi\cdot x+(1-\xi)\cdot elo_k^{cl}(x,p)=x+(1-\xi)\cdot (elo_k^{cl}(x,p)-x)
\end{equation}
moves ratings towards their images with respect to the classical Elo map. Hence, $x\in\mathcal{R}$ is a fixed point of $\varphi_\xi$ if and only if it is a self-justifying rating with respect to $p$ and $k$. In particular, the fixed points of the maps (\ref{phidef}) do not depend on $\xi$.
\end{Rem}

\begin{Thm}[Existence and uniqueness]\label{hauptthm}
Let $k>0$ and $p\in\mathcal{P}$.
\begin{enumerate}[(i)]
\item{Then there exists one and only one self-justifying rating $x\in\mathcal{R}$ with respect to $p$ and $k$.}\label{hauptthm_i}
\item{The self-justifying rating with respect to $p$ and $k$ is the unique fixed point of the maps (\ref{phidef}) for all $0\leq\xi<1$.}\label{hauptthm_ii}
\item{Let $\frac{G}{G+1}\leq\xi<1$ with $G:=\frac{k\cdot(n-1)}{4}\cdot \max_{i,j\in n}{(p_{ij}+p_{ji})}\geq0$. Then the map $\varphi_\xi$ in (\ref{phidef}) is a contraction with contraction factor $\xi$, i.e. we have
\begin{equation}\label{eqcontr}
\norm*{\varphi_\xi(z)-\varphi_\xi(y)}_1\leq\xi\cdot\norm*{ z-y}_1
\end{equation}
for all $y,z\in\mathcal{R}$.}\label{hauptthm_iii}
\item{Let $\frac{G}{G+1}\leq\xi<1$ and $y\in\mathcal{R}$. Then the sequence $\left(\varphi_\xi^l(y)\right)_{l\in\mathbb{N}}$ converges to the self-justifying rating $x$ with respect to $p$ and $k$. For $l\in\mathbb{N}$, we have
\begin{align}
\norm*{\varphi_\xi^l(y)-x}_1&\leq\xi^l\cdot\norm*{y-x}_1\nonumber\\
\norm*{\varphi_\xi^l(y)-\varphi_\xi^{l+1}(y)}_1&\leq\xi^l\cdot\norm*{y-\varphi_\xi(y)}_1\label{convphifolge}\\
\norm*{\varphi_\xi^l(y)-elo_k^{cl}(\varphi_\xi^{l}(y),p)}_1&\leq\xi^l\cdot\norm*{y-elo_k^{cl}(y,p)}_1\label{convelofolge}
\text{.}
\end{align}}\label{hauptthm_iv}
\end{enumerate}
\end{Thm}

\begin{Bew}
It will be enough to prove part (\ref{hauptthm_iii}): since $\mathcal{R}$ is a complete metric space (with respect to the metric induced by $\norm*{.}_1$), we can apply the Banach fixed point theorem to the contraction $\varphi_\xi$, and we obtain a unique fixed point. Together with Remark \ref{Remphi}, this yields (\ref{hauptthm_i}) and (\ref{hauptthm_ii}). Most of part (\ref{hauptthm_iv}) is an immediate formal consequence of (\ref{eqcontr}). Inequality (\ref{convelofolge}) follows from (\ref{convphifolge}) via
\begin{equation*}
\varphi_\xi^l(y)-\varphi_\xi^{l+1}(y)
=\varphi_\xi^l(y)-\xi\cdot\varphi_\xi^l(y)-(1-\xi)\cdot elo_k^{cl}(\varphi_\xi^l(y),p)
=(1-\xi)\cdot\left(\varphi_\xi^l(y)-elo_k^{cl}(\varphi_\xi^l(y),p)\right)\text{.}
\end{equation*}

To prove part (\ref{hauptthm_iii}), we have to extend the definition of $\varphi_\xi$. For $y\in\mathcal{R}$, we have
\begin{align}
\varphi_\xi(y)
&=\xi\cdot y+(1-\xi)\cdot elo_k^{cl}(y,p)\nonumber\\
&=\xi\cdot y+(1-\xi)\cdot k\cdot\left(\sum\limits_{j\in n}{\left(p_{ij}-\frac{p_{ij}+p_{ji}}{1+exp(y_j-y_i)}\right)}\right)_{i\in n}\text{.}\label{eqphidet}
\end{align}
We use expression (\ref{eqphidet}) to define a map $\varphi_\xi:\mathbb{R}^n\rightarrow\mathbb{R}^n$. We will show in two steps that (\ref{eqcontr}) holds for all $y,z\in\mathbb{R}^n$.

Step 1: We start by considering elements $y,z\in\mathbb{R}^n$ that differ by only one component (i.e. $\abs*{\{i\in n\ |\ y_i\neq z_i\}}=1$). By renumbering the index set, we can assume that $y_{n-1}<z_{n-1}$, whereas $y_i=z_i$ for $i\in n-1$ (use symmetry of definitions of $\varphi_\xi$ and $G$).

If we set
\begin{equation*}
T_i:\mathbb{R}\rightarrow\mathbb{R},\ t\mapsto-\frac{(1-\xi)\cdot k\cdot(p_{i(n-1)}+p_{(n-1)i})}{1+exp(t-y_i)}
\end{equation*}
for $i\in n-1$, then we see from (\ref{eqphidet}) and (\ref{expconst}) that 
\begin{align*}
(\varphi_\xi(z))_i-(\varphi_\xi(y))_i=\begin{cases}T_i(z_{n-1})-T_i(y_{n-1})&\text{if }i\in n-1\\
\xi\cdot(z_{n-1}-y_{n-1})-\sum\limits_{j\in n-1}{\left(T_j(z_{n-1})-T_j(y_{n-1})\right)}&\text{if }i=n-1\end{cases}
\end{align*}
for $i\in n$. Obviously, for all $i\in n-1$, we have $T_i(z_{n-1})-T_i(y_{n-1})\geq0$. Since the derivative of the map $T:\mathbb{R}\rightarrow\mathbb{R},\ t\mapsto\frac{1}{1+exp(t)}$ is bounded by $\max_{t\in\mathbb{R}}{\abs*{T'(t)}}=\frac{1}{4}$, the mean value theorem yields
\begin{align*}
&\sum\limits_{j\in n-1}{\left(T_j(z_{n-1})-T_j(y_{n-1})\right)}\\
=&\sum\limits_{j\in n-1}{(1-\xi)\cdot k\cdot(p_{j(n-1)}+p_{(n-1)j})\cdot\left(\frac{1}{1+exp(y_{n-1}-y_j)}-\frac{1}{1+exp(z_{n-1}-y_j)}\right)}\\
\leq&\sum\limits_{j\in n-1}{(1-\xi)\cdot k\cdot(p_{j(n-1)}+p_{(n-1)j})\cdot\frac{z_{n-1}-y_{n-1}}{4}}\\
\leq&(n-1)\cdot(1-\xi)\cdot k\cdot \max\limits_{i,j\in n}{(p_{ij}+p_{ji})}\cdot\frac{z_{n-1}-y_{n-1}}{4}\\
=&(1-\xi)\cdot G\cdot(z_{n-1}-y_{n-1})\\
\leq&\xi\cdot(z_{n-1}-y_{n-1})\text{,}
\end{align*}
where the last inequality follows from the assumption: $\frac{G}{G+1}\leq\xi$ implies $G\leq(G+1)\cdot\xi$ and $G\cdot(1-\xi)\leq\xi$. Now, we see that equality holds in (\ref{eqcontr}):
\begin{align*}
\norm*{\varphi_\xi(z)-\varphi_\xi(y)}_1
&=\sum\limits_{i\in n}{\abs*{\left(\varphi_\xi(z)\right)_i-\left(\varphi_\xi(y)\right)_i}}\\
&=\sum\limits_{i\in n}{\left(\left(\varphi_\xi(z)\right)_i-\left(\varphi_\xi(y)\right)_i\right)}&&\text{(all summands non-negative)}\\
&=\xi\cdot(z_{n-1}-y_{n-1})&&\text{(}T_i\text{-terms cancel out)}\\
&=\xi\cdot\norm*{z-y}_1\text{.}
\end{align*}

Step 2: We conclude by proving (\ref{eqcontr}) for general $y,z\in\mathbb{R}^n$. Define $y^l\in\mathbb{R}^n$ by $y^l_i:=\begin{cases}y_i\ &\text{if }l\leq i\\z_i\ &\text{if }l>i\end{cases}$ for $l\in(n+1),i\in n$. Then
\begin{align*}
\norm*{\varphi_\xi(z)-\varphi_\xi(y)}_1
&=\norm*{\varphi_\xi(y^n)-\varphi_\xi(y^0)}_1\\
&\leq\sum\limits_{l\in n}\norm*{\varphi_\xi(y^{l+1})-\varphi_\xi(y^l)}_1&&\text{(triangle inequality)}\\
&\leq\sum\limits_{l\in n}\xi\cdot\norm*{ y^{l+1}-y^l}_1&&\text{(by step 1)}\\
&=\xi\cdot\sum\limits_{l\in n}\abs*{y^{l+1}_l-y^l_l}\\
&=\xi\cdot\sum\limits_{l\in n}\abs*{z_l-y_l}\\
&=\xi\cdot\norm*{z-y}_1\text{.}
\end{align*}
\hfill$\Box$
\end{Bew}

\begin{Def}[Self-justifying Elo map]
Let $k>0$. The \emph{self-justifying Elo map} $elo_k:\mathcal{P}\rightarrow\mathcal{R}$ with dynamising parameter $k$ maps $p\in\mathcal{P}$ to the unique self-justifying rating with respect to $p$ and $k$ guaranteed by theorem \ref{hauptthm} (\ref{hauptthm_i}).
\end{Def}

\begin{App}[Self-justifying Elo]\label{Appsj}
We resume the setting of the classical application \ref{Appcl}, i.e. we have chosen $k>0$, an average rating $\mu\in\mathbb{R}$, and a deviation factor $\sigma>0$. Again, we consider the (finite) sequence $p^0,p^1,...,p^{m-1}\in\mathcal{P}$ of competition results that correspond to consecutive time periods.

We have two approaches in mind. In a time-limited context (e.g. a single tournament or a ranking with respect to a single year), all results should be equally weighted, and we set $q^l:=\sum_{i=0}^l{p^i}$ for $l\in m$. In a long-term scenario, it might be better to use weights favouring more current results. If the time periods are of equal duration, this can be implemented by choosing a suitable factor $0<f<1$ and letting $q^l:=\sum_{i=0}^l{f^{l-i}\cdot p^i}$ for $l\in m$.

In either case, for $l\in m$, the self-justifying Elo rating at the end of period $l$ is $elo_k(q^l)\in\mathcal{R}$. $elo_k(q^{m-1})$ is the current self-justifying Elo rating.

As in the classical case, we transform the result at the end, getting the published rating $(\mu)_{i\in n}+\sigma\cdot elo_k(q^{m-1})$.
\end{App}

\begin{Rem}[Initial ratings]
Elo describes several methods to generate intital ratings for new players. He introduces a `method of successive approximations' for a group of unrated players with existing results, and formulas for initial ratings resulting from special `rating tournaments' (\cite{ELO}, sections 3.3, 3.4, 1.7). We suggest not to use a differing procedure to provide initial ratings. Instead, new players should simply be added to the players pool and application \ref{Appsj} will do its job very naturally. However, it makes sense to define a required number of competition results against rated players that a new player must reach before he or she (and the corresponding results) are included in the calculation.

Side remark: The `method of successive approximations' comes closest to our self-justifying approach, since it aims at the convergence of a particular sequence of ratings to a fixed point. Nevertheless, it behaves much worse. For example, the definition of the sequence requires extra assumptions (no full or zero scores), and convergence cannot be guaranteed (not even in a weaker sense, only considering differences). Furthermore, it can happen that a win against a low-rated player leads to a lower rating. Unsurprisingly, Elo considers the procedure rather sceptically and restricts its field of application.
\end{Rem}

\begin{Exa}[Overcompensation]\label{Exaoco}
Let $k=1$. In a scenario with two players, consider the competition result $p\in\mathcal{P}$ defined by $p_{01}=55$ and $p_{10}=45$. A rating of $\ln{\left(\sqrt{\frac{55}{45}}\right)}\approx0.10$ for player $0$ would reflect this result perfectly, meaning that the adjustments of the classical Elo \ref{Appcl} would leave such a rating unchanged. Starting from the initial rating $0\in\mathcal{R}$, it is clear that player $0$'s rating should slightly increase. Classical Elo, however, changes it to $5$, indicating a significant overcompensation.

How can such an inaccuracy happen? Of course, the large number of points that are evaluated at the same time contributes to the result (as does the choice of the parameter $k$). But the essential reason is the fact that the classical Elo process only considers performance expectations with respect to the previous rating. So the performance expectations with respect to the adjusted rating are completely ignored -- they are far from corresponding to the actual competition result.

If we confront the self-justifying rating with this example, it will return a rating near the `perfect' one above.
\end{Exa}

\begin{Lem}[Estimation of precision]\label{Lemprec}
Let $k>0$ and $p\in\mathcal{P}$. Let $x\in\mathcal{R}$. Then
\begin{equation*}
\norm*{x-elo_k(p)}_1\leq\norm*{x-elo_k^{cl}(x,p)}_1
\end{equation*}
\end{Lem}

\begin{Bew}
Let $\xi$ be as in theorem \ref{hauptthm} (\ref{hauptthm_iii}). Then the map $\varphi_\xi$ defined in (\ref{phidef}) is a contraction with contraction factor $\xi$. We find that
\begin{align*}
&\norm*{x-elo_k(p)}_1\\
\leq&\norm*{x-\varphi_\xi(x)}_1+\norm*{\varphi_\xi(x)-elo_k(p)}_1&&\text{(triangle inequality)}\\
=&\norm*{x-\xi\cdot x-(1-\xi)\cdot elo_k^{cl}(x,p)}_1+\norm*{\varphi_\xi(x)-\varphi_\xi(elo_k(p))}_1\\
\leq&(1-\xi)\cdot\norm*{x-elo_k^{cl}(x,p)}_1+\xi\cdot\norm*{x-elo_k(p)}_1&&\text{(contraction)}
\text{.}
\end{align*}
Subtraction yields $(1-\xi)\cdot\norm*{x-elo_k(p)}_1\leq(1-\xi)\cdot\norm*{x-elo_k^{cl}(x,p)}_1$ and the assertion of the lemma.
\hfill$\Box$
\end{Bew}

\section{Properties}

\begin{Prp}[Continuity]\label{Prpcont}
The map
\begin{equation*}
\mathcal{P}\times\mathbb{R}_{>0}\rightarrow\mathcal{R}:(p,k)\rightarrow elo_k(p)
\end{equation*}
is continuous. For $p,q\in\mathcal{P}$ and $k>0$, we have
\begin{align}
\norm*{elo_k(p)-elo_k(q)}_1&\leq 2\cdot k\cdot\norm*{p-q}_1\label{cont1}
\text{.}
\end{align}
\end{Prp}

\begin{Bew}
We start by proving (\ref{cont1}). An argument similar to the one in step 2 of the proof of theorem \ref{hauptthm} shows that we can assume that $p$ and $q$ only differ by one component, say $p_{01}<q_{01}$, whereas all other components coincide. Then by definition the classical Elo map,
\begin{align*}
\norm*{elo_k^{cl}(x,p)-elo_k^{cl}(x,q)}_1&=\sum\limits_{i=0}^1{\abs*{\left(elo_k^{cl}(x,p)\right)_i-\left(elo_k^{cl}(x,q)\right)_i}}&&\text{(others unaffected)}\\
&=2\cdot\abs*{\left(elo_k^{cl}(x,p)\right)_0-\left(elo_k^{cl}(x,q)\right)_0}&&\text{(by remark \ref{Remsym})}\\
&=2\cdot k\cdot\abs*{(p_{01}-q_{01})\cdot\left(1-\frac{1}{1+exp(x_j-x_i)}\right)}\\
&\leq2\cdot k\cdot\abs*{p_{01}-q_{01}}\\
&=2\cdot k\cdot\norm*{p-q}_1
\end{align*}
for $x\in\mathcal{R}$. The required inequality
\begin{align*}
\norm*{elo_k(p)-elo_k(q)}_1
&\leq\norm*{elo_k(p)-elo_k^{cl}(elo_k(p),q)}_1&&\text{(by lemma \ref{Lemprec})}\\
&=\norm*{elo_k^{cl}(elo_k(p),p)-elo_k^{cl}(elo_k(p),q)}_1&&\text{(definition of }elo_k\text{)}\\
&\leq 2\cdot k\cdot\norm*{p-q}_1
\end{align*}
follows.

The asserted continuity follows from the factorization
\begin{equation*}
\mathcal{P}\times\mathbb{R}_{>0}\xrightarrow{(p,k)\mapsto k\cdot p}\mathcal{P}\xrightarrow{elo_1}\mathcal{R}
\end{equation*}
(use remark \ref{Remsuppk}), where continuity of the first map is obvious, and continuity of the second map follows from (\ref{cont1}).
\hfill$\Box$
\end{Bew}

\begin{Rem}[Self-correction]\label{Remsel}
If player $i$ wins one point against player $j$, then the classical application \ref{Appcl} of Elo increases player $i$'s rating by at most $k$ and decreases player $j$'s rating by the same value, hence the total change of the rating is at most $2\cdot k$. Inequality (\ref{cont1}) shows that the same holds true for the self-justifying rating.

Here, however, not only the ratings of the involved players $i$ and $j$ may change. Assume that a third player has had several games against player $i$ before. Taking into account player $i$'s modified rating, the third player's performance in these games has to be regarded more generously -- so his or her rating should also be increased. The self-justifying approach follows this idea in a natural way.

The classical Elo encourages players to choose their opponents selectively. Assume that a player is underrated. Then others will avoid to compete with him or her, because they cannot expect a fair adjustment of the rating.

The phenomenon seen above makes the self-justifying Elo behave better: when the underrated player regains a realistic rating afterwards, the unfair adjustment is corrected automatically.

Side remark: Elo also describes his system as self-correcting (\cite{ELO}, 1.67). What he means is, for example, that the competition results of a player whose rating is lower than the `correct' one will be compared with too low performance expectations. Hence the rating will tend to increase.
\end{Rem}

\begin{Prp}[Monotonicity]\label{Prpmon}
Let $k>0$. Let $p,q\in\mathcal{P}$ be competition results with $q_{i'j'}>p_{i'j'}$ for one pair $(i',j')\in n\times n$ and $q_{ij}=p_{ij}$ for all $(i,j)\in n\times n\setminus\{(i',j')\}$. Then $\left(elo_k(q)\right)_{i'}>\left(elo_k(p)\right)_{i'}$ and $\left(elo_k(q)\right)_{j'}<\left(elo_k(p)\right)_{j'}$.
\end{Prp}

\begin{Bew}
By theorem \ref{hauptthm} (\ref{hauptthm_iii}), we can choose $0\leq\xi<1$ in a way such that the map
\begin{equation*}
\varphi_\xi:\mathcal{R}\rightarrow\mathcal{R},\ x\mapsto x+(1-\xi)\cdot (elo_k^{cl}(x,q)-x)
\end{equation*}
is a contraction. Let $x:=elo_k(p)$ and $y:=elo_k(q)$. Then
\begin{align*}
\varphi_\xi(x)=&x+(1-\xi)\cdot (elo_k^{cl}(x,q)-x)\\
=&x+(1-\xi)\cdot (elo_k^{cl}(x,q)-elo_k^{cl}(x,p))&&\text{(definition of }elo_k\text{)}\\
=&x+\left(\begin{cases}C&\text{if }i=i'\\-C&\text{if }i=j'\\0&\text{otherwise}\end{cases}\right)_{i\in n}
\text{,}
\end{align*}
where $C:=(1-\xi)\cdot k\cdot(q_{i'j'}-p_{i'j'})\cdot\left(1-\frac{1}{1+exp(x_{j'}-x_{i'})}\right)>0$.

By increasing $\xi$ if necessary, we can assume that $y_{i'}\notin(x_{i'},x_{i'}+C)$ and that $y_{j'}\notin(x_{j'}-C,x_{j'})$. Since $\varphi_\xi$ is a contraction, $y$ is it fixed point, and $x$ is not its fixed point (as seen above), it follows that 
\begin{align*}
0&<\norm*{x-y}_1-\norm*{\varphi_\xi(x)-y}_1\\
&=\left(\abs*{x_{i'}-y_{i'}}-\abs*{x_{i'}+C-y_{i'}}\right)+\left(\abs*{x_{j'}-y_{j'}}-\abs*{x_{j'}-C-y_{j'}}\right)\\
&=\left(\begin{cases}C&\text{if }y_{i'}\geq x_{i'}+C\\-C&\text{if }y_{i'}\leq x_{i'}\end{cases}\right)
+\left(\begin{cases}C&\text{if }y_{j'}\leq x_{j'}-C\\-C&\text{if }y_{j'}\geq x_{j'}\end{cases}\right)
\text{.}
\end{align*}
Hence both summands are positive and the assertion follows.
\hfill$\Box$
\end{Bew}

\begin{Exa}[Lack of monotonicity]\label{Examon}
Let $k=1$, and consider a scenario with two players and two time periods. In period $0$, player $0$ wins one point ($p_{01}^0=1$, $p_{10}^0=0$), whereas player $1$ wins three points in period $1$ ($p_{01}^1=0$, $p_{10}^1=3$). According to application \ref{Appcl}, player $0$'s classical rating at the end is approximately $-1.69$. If the result of period $0$ is ignored, however, player $0$'s classical rating improves to $-1.5$.
\end{Exa}

\begin{Prp}[Behaviour for small $k$]\label{Prpsmall}
Let $p\in\mathcal{P}$ and $l>0$. Then
\begin{equation*}
\lim\limits_{k\downarrow 0}{elo_k(p)}=\lim\limits_{k\downarrow 0}{elo_l(k\cdot p)}=0
\end{equation*}
and
\begin{equation*}
\lim\limits_{k\downarrow 0}{\frac{elo_k(p)}{k}}=\lim\limits_{k\downarrow 0}{\frac{elo_l(k\cdot p)}{k\cdot l}}=\left(\sum\limits_{j\in n}{\frac{p_{ij}-p_{ji}}{2}}\right)_{i\in n}\text{.}
\end{equation*}
\end{Prp}

\begin{Bew}
For all $x\in\mathbb{R}$ and $k>0$, we have
\begin{align*}
\norm*{elo_k^{cl}(x,p)}_1
&=\norm*{k\cdot\left(\sum\limits_{j\in n}{\left(p_{ij}-\frac{p_{ij}+p_{ji}}{1+exp(x_j-x_i)}\right)}\right)_{i\in n}}_1\\
&\leq k\cdot\sum\limits_{i,j\in n}{\left(\abs*{p_{ij}}+\abs*{p_{ij}+p_{ji}}\right)}
\text{,}
\end{align*}
where the second factor only depends on $p$. Since $elo_k(p)=elo_k^{cl}(elo_k(p),p)$ by definition of $elo_k$, the first assertion results immediately (use $elo_l(k\cdot p)=elo_{k\cdot l}(p)$). Using this, we get
\begin{align*}
\lim\limits_{k\downarrow 0}{\frac{elo_k(p)}{k}}
&=\lim\limits_{k\downarrow 0}{\frac{elo_k^{cl}(elo_k(p),p)}{k}}\\
&=\lim\limits_{k\downarrow 0}{\left(\sum\limits_{j\in n}{\left(p_{ij}-\frac{p_{ij}+p_{ji}}{1+exp(\left(elo_k(p)\right)_j-\left(elo_k(p)\right)_i)}\right)}\right)_{i\in n}}\\
&=\left(\sum\limits_{j\in n}{\frac{p_{ij}-p_{ji}}{2}}\right)_{i\in n}
\text{.}
\end{align*}
\hfill$\Box$
\end{Bew}

\begin{Def}[Connectivity]
Let $p\in\mathcal{P}$. The \emph{result graph} corresponding to $p$ is the directed graph with vertice set $n$ and an edge from $i\in n$ to $j\in n$ whenever $p_{ij}>0$.

$p$ is said to be \emph{weakly connected} if its result graph is weakly connected, i.e. when there exists a not necessarily directed path from $i$ to $j$ for every $i,j\in n$. The \emph{connected components} of $p$ are the vertice sets $t\subseteq n$ of the maximal weakly connected subgraphs of the result graph. The \emph{restriction} of $p$ to its connected component $t$ is the competition result $p^t:=\left(\begin{cases}p_{i,j}&\text{if }i,j\in t\\0&\text{otherwise}\end{cases}\right)_{(i,j)\in n\times n}\in\mathcal{P}$.

$p$ is said to be \emph{strongly connected} if its result graph is strongly connected, i.e. when there exists a directed path from $i$ to $j$ for every $i,j\in n$. A connected component $t\subseteq n$ of $p$ is strongly connected if the subgraph of $p$ corresponding to the vertice set $t$ is strongly connected.
\end{Def}

\begin{Rem}[Connected components]\label{Remcon}
Let $p\in\mathcal{P}$, $k>0$. Let $(t_i)_{i\in m}$ be the connected components with restrictions $p^{t_i}$ of $p$. Then it is obvious that $p=\sum_{i\in m}{p^{t_i}}$ and $elo_k^{cl}(.,p)=\sum_{i\in m}{elo_k^{cl}(.,p^{t_i})}$. It is clear by definition that $elo_k(p^{t_i})$ can only have non-zero components within $t_i$. Hence 
\begin{align*}
elo_k^{cl}\left(\sum\limits_{i\in m}{elo_k(p^{t_i})},p\right)
&=\sum\limits_{j\in m}{elo_k^{cl}\left(\sum\limits_{i\in m}{elo_k(p^{t_i})},p^{t_j}\right)}\\
&=\sum\limits_{j\in m}{elo_k^{cl}(elo_k(p^{t_j}),p^{t_j})}\\
&=\sum\limits_{i\in m}{elo_k(p^{t_i})}&&\text{(definition of }elo_k\text{)}
\text{,}
\end{align*}
which implies $elo_k(p)=\sum_{i\in m}{elo_k(p^{t_i})}$ by theorem \ref{hauptthm} (\ref{hauptthm_i}).

Bearing in mind remark \ref{Rememb}, this allows us to assume weak connectivity in many situations.
\end{Rem}

\begin{Thm}[Behaviour for large $k$]\label{Thmbigk}
Let $p\in\mathcal{P}$. The following are equivalent.
\begin{enumerate}[(i)]
\item{The limit $\lim\limits_{k\rightarrow\infty}{elo_k(p)}$ exists in $\mathcal{R}$.}\label{large_b}
\item{For all $k>0$, the limit $\lim\limits_{l\rightarrow\infty}{elo_k(l\cdot p)}$ exists in $\mathcal{R}$.}\label{large_b1}
\item{The set $\{elo_k(p)\ |\ k>0\}\subseteq\mathcal{R}$ is bounded.}\label{large_e}
\item{There exists an element $x\in\mathcal{R}$ with $elo_1^{cl}(x,p)=0$.}\label{large_d}
\item{All connected components of $p$ are strongly connected.}\label{large_c}
\end{enumerate}

The following hold true.
\begin{enumerate}[(i)]
\setcounter{enumi}{5}
\item{If the equivalent conditions above are satisfied, then the limits in (\ref{large_b}) and (\ref{large_b1}) coincide and satisfy condition (\ref{large_d}).}\label{large_ext1}
\item{If $p$ is weakly connected, then there can exist at most one element $x\in\mathcal{R}$ satisfying condition (\ref{large_d}).}\label{large_ext2}
\end{enumerate}
\end{Thm}

\begin{Bew}
We start by proving extra assertion (\ref{large_ext2}). Assume that $p$ is weakly connected, and consider $x,y\in\mathcal{R}$ satisfying condition (\ref{large_d}), i.e. $elo_1^{cl}(x,p)=0=elo_1^{cl}(y,p)$. Let $t\subseteq n$ be the set of all $i\in n$ where $y_i-x_i$ is maximal. Then
\begin{align*}
0&=\sum\limits_{i\in t}{\left(\left(elo_1^{cl}(x,p)\right)_i-\left(elo_1^{cl}(y,p)\right)_i\right)}\\
&=\sum\limits_{i\in t}{\left(\sum\limits_{j\in n}{\left(p_{ij}-\frac{p_{ij}+p_{ji}}{1+exp(x_j-x_i)}\right)}-\sum\limits_{j\in n}{\left(p_{ij}-\frac{p_{ij}+p_{ji}}{1+exp(y_j-y_i)}\right)}\right)}\\
&=\sum\limits_{\substack{i\in t\\j\in n}}{(p_{ij}+p_{ji})\cdot\left(\frac{1}{1+exp(y_j-y_i)}-\frac{1}{1+exp(x_j-x_i)}\right)}
\text{.}
\end{align*}
By construction of $t$, every summand is non-negative ($y_i-x_i\geq y_j-x_j$ implies $y_j-y_i\leq x_j-x_i$), and hence it is equal to $0$. It follows that for every $i\in t,j\in n\setminus t$, we have $p_{ij}+p_{ji}=0$. Thus $t$ and $n\setminus t$ are not connected. Since $p$ is weakly connected and $t$ is non-empty, it follows that $t=n$ and thus $x=y$.

(\ref{large_b}) $\iff$ (\ref{large_b1}): By definition, $elo_k(l\cdot p)=elo_{k\cdot l}(p)$.

(\ref{large_b}) $\implies$ (\ref{large_d}): Let $x:=\lim_{k\rightarrow\infty}{elo_k(p)}$. By definition, we have 
\begin{equation*}
elo_k(p)=elo_k^{cl}(elo_k(p),p)=k\cdot elo_1^{cl}(elo_k(p),p)
\end{equation*}
for $k>0$. Because of the assumed convergence, it follows that $elo_1^{cl}(elo_k(p),p)$ converges to $0$. By continuity of $elo_1^{cl}$, however, this expression also converges to $elo_1^{cl}(x,p)$.

This also finishes the proof of extra assertion (\ref{large_ext1}).

(\ref{large_d}) $\implies$ (\ref{large_c}): Consider an arbitrary element $m\in n$ and let $t\subseteq n$ be the set of $i\in n$ that can be reached by a directed path beginning in $m$. We find that
\begin{align*}
0&=\sum\limits_{i\in t}{\left(elo_1^{cl}(x,p)\right)_i}\\
&=\sum\limits_{i\in t}{\sum\limits_{j\in n}{\left(p_{ij}-\frac{p_{ij}+p_{ji}}{1+exp(x_j-x_i)}\right)}}\\
&=\sum\limits_{\substack{i\in t\\j\in n\setminus t}}{\left(p_{ij}-\frac{p_{ij}+p_{ji}}{1+exp(x_j-x_i)}\right)}&&\text{(by remark \ref{Remsym})}\\
&=-\sum\limits_{\substack{i\in t\\j\in n\setminus t}}{\left(\frac{p_{ji}}{1+exp(x_j-x_i)}\right)}&&\text{(by construction of }t\text{)}
\text{.}
\end{align*}
It follows that $p_{ji}=0$ for all $i\in t,j\in n\setminus t$, showing that $t$ is a connected component of $p$. Since $m$ was chosen arbitrarily, the connected components of $p$ are strongly connected.

(\ref{large_c}) $\implies$ (\ref{large_e}): By remark \ref{Remcon}, it is enough to prove boundedness for the connected components. We will assume that $p$ is strongly connected.

For every proper non-empty subset $t\subset n$, we choose $C^t>0$ such that
\begin{equation}\label{hilf_bound_1}
\sum\limits_{\substack{i\in t\\j\in n\setminus t}}{\left(p_{ij}-\frac{p_{ij}+p_{ji}}{1+exp(-C)}\right)}<0
\end{equation}
for all $C\geq C^t$. This is possible because of the strong connectedness: there must exist $i\in t,j\in n\setminus t$ with $p_{ji}>0$, and hence the limit of the decreasing expression (\ref{hilf_bound_1}) for $C\rightarrow\infty$ is at most $-p_{ji}$. We set $C:=\max_{\emptyset\subsetneq t\subsetneq n}{C^t}$.

Now let $k>0$ and set $y:=elo_k(p)$. We will show that for $j',i'\in n$ with components of consecutive size (meaning that $y_{j'}<y_{i'}$ and that no other component of $y$ takes a value between them), we have $y_{i'}-y_{j'}<C$. The assertion follows immediately (using that $y\in\mathcal{R}$).

Let $t:=\{i\in n\ |\ y_i\geq y_{i'}\}\subset n$. Then
\begin{align*}
0&\leq\sum\limits_{i\in t}{y_i}\\
&=\sum\limits_{i\in t}{\left(elo_k^{cl}(y,p)\right)_{i}}&&\text{(definition of }elo_k\text{)}\\
&=\sum\limits_{i\in t}{k\cdot\sum\limits_{j\in n}{\left(p_{ij}-\frac{p_{ij}+p_{ji}}{1+exp(y_j-y_i)}\right)}}\\
&=k\cdot\sum\limits_{\substack{i\in t\\j\in n\setminus t}}{\left(p_{ij}-\frac{p_{ij}+p_{ji}}{1+exp(y_j-y_i)}\right)}&&\text{(by remark \ref{Remsym})}\\
&\leq k\cdot\sum\limits_{\substack{i\in t\\j\in n\setminus t}}{\left(p_{ij}-\frac{p_{ij}+p_{ji}}{1+exp(y_{j'}-y_{i'})}\right)}&&\text{(by construction of }t\text{)}
\text{.}
\end{align*}
It follows from (\ref{hilf_bound_1}) that $y_{i'}-y_{j'}<C^t\leq C$ as required.

(\ref{large_e}) $\implies$ (\ref{large_b}): $elo_k(p)$ converges / is bounded if and only if $elo_k(p^t)$ converges / is bounded for all connected components $t$ of $p$. Hence we can assume that $p$ is weakly connected.

Consider a sequence $(k_i)_{i\in\mathbb{N}}\in\mathbb{R}_{>0}^\mathbb{N}$ with $\lim_{i\rightarrow\infty}{k_i}=\infty$ for which $(elo_{k_i}(p))_{i\in\mathbb{N}}$ converges to some element $x\in\mathcal{R}$. By Bolzano-Weierstraß, $elo_k(p)$ converges (for $k\rightarrow\infty$) if and only if the limits $x$ of all such sequences coincide.

Since $\{elo_k(p)\ |\ k>0\}$ is bounded, we have
\begin{align*}
elo_1^{cl}(x,p)=\lim\limits_{i\rightarrow\infty}{elo_1^{cl}(elo_{k_i}(p),p)}=\lim\limits_{i\rightarrow\infty}{\frac{1}{k_i}\cdot elo_{k_i}^{cl}(elo_{k_i}(p),p)}=\lim\limits_{i\rightarrow\infty}{\frac{1}{k_i}\cdot elo_{k_i}(p)}=0
\text{.}
\end{align*}

By extra assertion (\ref{large_ext2}) proven above, there exists at most one point $x\in\mathcal{R}$ satisfying this condition (use weak connectedness), which concludes the proof of the theorem.
\hfill$\Box$
\end{Bew}

\begin{Kor}[Classical convergence]\label{Korconv}
Let $p\in\mathcal{P}$, $k>0$. Define a sequence of competition results by $p^l:=p$ for $l\in\mathbb{N}$. Let $x^l\in\mathcal{R}$ be the classical Elo rating after period $l$ as defined in application \ref{Appcl}, i.e. $x^0=0$ and $x^l=x^{l-1}+elo_k^{cl}(x^{l-1},p)$ for $l\in\mathbb{N}\setminus\{0\}$. If the sequence $(x^l)_{l\in\mathbb{N}}$ of classical Elo ratings converges to a rating $x\in\mathcal{R}$, then $x=\lim_{l\rightarrow\infty}elo_l(p)=\lim_{l\rightarrow\infty}elo_k(l\cdot p)$.
\end{Kor}

\begin{Bew}
By remark \ref{Remcon}, we can assume weak connectedness of $p$.

Convergence implies $elo_k^{cl}(x,p)=\lim_{l\rightarrow\infty}{elo_k^{cl}(x^l,p)}=0$ by continuity of $elo_k^{cl}$ and by definition of the sequence. Theorem \ref{Thmbigk} (\ref{large_d}) $\implies$ (\ref{large_b}), (\ref{large_ext1}), (\ref{large_ext2}) concludes the proof.
\hfill$\Box$
\end{Bew}

\begin{Exa}[Lack of convergence]\label{Exadiv}
Let $k=1$. We discuss another example with two players. A competition result $p\in\mathcal{P}$ given by $p_{01}=3$ and $p_{10}=2$ is continuously repeated, i.e. $p^l:=p$ for $l\in\mathbb{N}$. By definition of classical Elo \ref{Appcl}, the map $e:\mathbb{R}\rightarrow\mathbb{R},\ y\mapsto y+3-\frac{5}{1+exp(-2\cdot y)}$ determines the change of the $0$-th component of the rating within one time period. It is clear that $p$ is strongly connected, so we know from theorem \ref{Thmbigk} (\ref{large_c}) $\implies$ (\ref{large_d}), (\ref{large_ext2}) that there exists a unique $x\in\mathcal{R}$ with $elo_1^{cl}(x,p)=0$. This corresponds to a unique fixed point $x_0=\ln{(\sqrt{\frac{3}{2}})}\approx0.20$ of $e$. However, the fixed point turns out to be repulsive, and the sequence of classical Elo ratings does not converge. Instead, there exist two additional attractive fixed points of $e\circ e$ near $1.05$ and $-0.40$, respectively, and the rating jumps between neighbourhoods of these points.
\end{Exa}

\section{Computation}\label{SctCalc}

\begin{Alg}[Self-justifying Elo]\label{alg2}
The following algorithm computes the self-justifying Elo.

Input: dynamising parameter $k>0$, competition result $p\in\mathcal{P}$, expected precision $\varepsilon>0$

Output: $x\in\mathcal{R}$ with $\norm*{x-elo_k(p)}_1\leq\varepsilon$

Technical parameter: continuity parameter $c\in\mathbb{N}$

\begin{enumerate}[(1)]
\item{Set $x:=0\in\mathcal{R}$. Set $\xi:=\frac{G}{G+1}$ with $G:=\frac{k\cdot(n-1)}{4}\cdot \max\limits_{i,j\in n}{(p_{ij}+p_{ji})}$. Set $d^{perm}:=\infty$, $w:=0$.}\label{alg2_}
\item{Set $u_{ij}:=k\cdot\left(p_{ij}-\frac{p_{ij}+p_{ji}}{1+exp(x_j-x_i)}\right)=-k\cdot\left(p_{ji}-\frac{p_{ij}+p_{ji}}{1+exp(x_i-x_j)}\right)$ for $i,j\in n$, $i<j$.}\label{alg2a}
\item{Set $e_i:=-\sum\limits_{j\in n,j<i}{u_{ji}}+\sum\limits_{j\in n,j>i}{u_{ij}}$ for $i\in n$.}\label{alg2b}
\item{Set $d:=\norm*{x-e}_1$. If $d\leq\varepsilon$, then return $x$.}\label{alg2c}
\item{If $d>\xi\cdot d^{perm}$, then set $\xi:=\sqrt{\xi}$, $w:=c$, and go to (\ref{alg2f}).}\label{alg2d}
\item{If $w>0$, then set $w:=w-1$. Otherwise, set $\xi:=\xi^2$.}\label{alg2e}
\item{If $d<d^{perm}$, then set $x^{perm}:=x$, $e^{perm}:=e$, $d^{perm}:=d$.}\label{alg2f}
\item{Set $x:=\xi\cdot x^{perm}+(1-\xi)\cdot e^{perm}$. Go to (\ref{alg2a}).}\label{alg2g}
\end{enumerate}
\end{Alg}

\begin{Bew}
After the initializations of step (\ref{alg2_}), the algorithm performs a loop consisting of steps (\ref{alg2a}) to (\ref{alg2g}). Steps (\ref{alg2a}) and (\ref{alg2b}) calculate $e=elo_k^{cl}(x,p)$ (use remark \ref{Remsym} to see that either formula in step (\ref{alg2a}) works; for computational reasons, it might be better to use the one with the negative argument of $exp$). Steps (\ref{alg2d}) and (\ref{alg2e}) adjust the variable $\xi$ depending on the change of $d=\norm*{x-e}_1$. Step (\ref{alg2g}) computes a new value of $x$, applying the function $\varphi_\xi$ defined in (\ref{phidef}).

Clearly, $x\in\mathcal{R}$ throughout the execution. The algorithm can only finish in step (\ref{alg2c}), where the condition and lemma \ref{Lemprec} guarantee the correctness of the result. We have to prove that the algorithm terminates eventually.

Step (\ref{alg2f}) generates a sequence of ratings $x^{perm}$ for which the values $d^{perm}=\norm*{x^{perm}-e^{perm}}_1=\norm*{x^{perm}-elo_k^{cl}(x^{perm},p)}_1$ decrease. Whenever $\xi=\frac{G}{G+1}$, it follows from inequality (\ref{convelofolge}) that $d\leq\xi\cdot d^{perm}$ in step (\ref{alg2d}) -- and thus the algorithm does not increase $\xi$ further. Hence $\frac{G}{G+1}$ is the maximal value that $\xi$ can take, and all loops in which $\xi$ is not increased in step (\ref{alg2d}) lead to a reduction of $d^{perm}$ of at least a factor of $\frac{G}{G+1}$. On the other hand, we have seen that $\xi$ can only be increased a finite number of times in sequence (until it reaches $\frac{G}{G+1}$ again). We deduce that the sequence of $d^{perm}$ converges to $0$, which implies that the algorithm terminates.
\hfill$\Box$
\end{Bew}

\begin{Rem}[Algorithmic approach]
The idea behind algorithm \ref{alg2} is simple. We know that repeated application of the function $\varphi_\xi$ defined in (\ref{phidef}) guarantees convergence for $\xi\geq\frac{G}{G+1}$. Lemma \ref{Lemprec} allows us to measure the precision of our intermediate results, and inequality (\ref{convelofolge}) determines our expectations concerning convergence speed.

Armed with these tools, we can simply try out smaller values of $\xi$ (hoping for some speed-up) -- we get an immediate feedback whether the choice is suitable or not. It is not obvious in what way our attempts can be exercised most efficiently. Algorithm \ref{alg2} uses a parameter $c$ to make sure that not too much time is wasted by failed attempts. It implements a period of continuity (no reduction of $\xi$) whenever $\xi$ does not meet the expectations.
\end{Rem}

\begin{Prp}[Execution duration]
The execution of algorithm \ref{alg2} takes at most
\begin{equation}\label{durationalg2}
\begin{cases}1&\text{if }2\cdot k\cdot\norm*{p}_1\leq\varepsilon\\
\left\lceil\left\lceil\left(\frac{k\cdot(n-1)}{4}\cdot \max\limits_{i,j\in n}{(p_{ij}+p_{ji})}+1\right)\cdot\ln{\left(\frac{2\cdot k\cdot\norm*{p}_1}{\varepsilon}\right)}\right\rceil\cdot\frac{c+2}{c+1}\right\rceil+1&\text{otherwise}\end{cases}
\end{equation}
loops (steps (\ref{alg2a}) to (\ref{alg2g})). Each loop includes at most $\frac{n^2}{2}+\frac{7\cdot n}{2}+4$ assignments of scalar variables.
\end{Prp}

\begin{Bew}
The maximal number of variable assignments in steps (\ref{alg2a}) to (\ref{alg2g}) are $\frac{n\cdot (n-1)}{2}$, $n$, $1$, $2$, $1$, $2\cdot n+1$, and $n$, respectively. The assignments in steps (\ref{alg2d}) and (\ref{alg2e}) cannot occur together. Summation yields the second assertion of the proposition.

The value of $d^{perm}$ after the first loop is $\norm*{0-elo_k^{cl}(0,p)}_1$. It follows from (\ref{absch_2kp}) that $d^{perm}\leq 2\cdot k\cdot\norm*{p}_1$. In the following, we assume that $2\cdot k\cdot\norm*{p}_1>\varepsilon$. Otherwise, the algorithm clearly terminates in the first loop. We count the first loop separately, it is represented by the `$+1$'-term in (\ref{durationalg2}).

We have seen in the proof of algorithm \ref{alg2} that the maximal value the variable $\xi$ can take is $\frac{G}{G+1}$, and that all loops in which $\xi$ is not increased lead to a reduction of $d^{perm}$ of at least a factor of $\frac{G}{G+1}$.

Next, consider the case in which $\xi$ is increased (say its previous value $\xi'$ is replaced by $\xi''$), but the resulting value $\xi''$ is still smaller than $\frac{G}{G+1}$. In this loop, $d^{perm}$ need not decrease. However, there has been a corresponding loop before in which $\xi''$ has been replaced by $\xi'$. This implies that $d^{perm}$ has decreased by at least a factor of $\xi''\leq\left(\frac{G}{G+1}\right)^2$. On average, both loops together guarantee the same reduction of $d^{perm}$ as do the regular loops without increase of $\xi$.

After $\left\lceil\ln{\left(\frac{2\cdot k\cdot\norm*{p}_1}{\varepsilon}\right)}\cdot\left(\ln{\left(\frac{G+1}{G}\right)}\right)^{-1}\right\rceil$ of the loops discussed before (not including the first one), the algorithm must have terminated. We can simplify this expression by using $\left(\ln{\left(\frac{G+1}{G}\right)}\right)^{-1}=\left(-\ln{\left(\frac{G}{G+1}\right)}\right)^{-1}=\left(-\ln{\left(1-\frac{1}{G+1}\right)}\right)^{-1}\leq\left(\frac{1}{G+1}\right)^{-1}=G+1$, where the inequality follows from the power series expansion $\ln(1-z)=-\sum_{i=1}^\infty{\frac{z^i}{i}}$ for $z\in(0,1)$.

We cannot guarantee any decrease of $d^{perm}$ for the remaining loops, i.e. those in which $\xi$ is increased from $\left(\frac{G}{G+1}\right)^2$ to $\frac{G}{G+1}$. The parameter $c$ makes sure that at most one in $c+2$ loops is of this type, possibly starting in loop $2$. So after a total number of $\left\lceil\left\lceil\ln{\left(\frac{2\cdot k\cdot\norm*{p}_1}{\varepsilon}\right)}\cdot(G+1)\right\rceil\cdot\frac{c+2}{c+1}\right\rceil+1$ loops, the required number of $d^{perm}$-reducing loops is reached.
\hfill$\Box$
\end{Bew}

\begin{Kor}[Asymptotic runtime]
Let $k>0$, $\varepsilon>0$, and $c\in\mathbb{N}$ be fixed. Assume that there exists an overall bound $C>0$ for the number of points for which a single player can compete, i.e. only competition results $p$ with $\sum_{j\in n}{\left(p_{ij}+p_{ji}\right)}\leq C$ for $i\in n$ need to be considered. Then the maximal number of assignments of scalar variables performed by algorithm \ref{alg2} for varying numbers $n$ of competing players is $\mathcal{O}\left(n^3\cdot\ln{\left(n\right)}\right)$.
\end{Kor}

\begin{Bew}
The main expression of (\ref{durationalg2}) is not larger than
\begin{equation*}
\left\lceil\left\lceil\left(\frac{k\cdot(n-1)}{4}\cdot C+1\right)\cdot\ln{\left(\frac{2\cdot k\cdot n\cdot C}{\varepsilon}\right)}\right\rceil\cdot\frac{c+2}{c+1}\right\rceil+1
\end{equation*}
(assume again that $2\cdot k\cdot\norm*{p}_1>\varepsilon$), and thus $\mathcal{O}\left(n\cdot\ln{(n)}\right)$. The maximal number of variable assignments per loop is $\mathcal{O}\left(n^2\right)$. The initial step (\ref{alg2_}) of the algorithm takes only $\mathcal{O}\left(n\right)$ variable assignments.
\hfill$\Box$
\end{Bew}

\end{document}